\documentclass[10pt]{article}
 
\usepackage{amsfonts}
\usepackage{amsmath}
\usepackage{longtable}
\usepackage{hyperref}
\usepackage{multirow}
\usepackage{booktabs}
\usepackage{paralist}
\usepackage{graphicx}
\usepackage[dvipsnames]{xcolor}

\newcommand{\halmos}{\hfill$\Box$}

\newcommand{\R}{\mathbb{R}}
\newcommand{\N}{\mathbb{N}}

\newcommand{\kbar}{\underline{k}}
\newcommand{\sigmin}{\sigma_{\min}}
\newcommand{\sigmax}{\sigma_{\max}}
\newcommand{\taumin}{\tau_{\min}}
\newcommand{\taumax}{\tau_{\max}}
\newcommand{\spg}{\mathrm{spg}}
\newcommand{\acc}{\mathrm{accel}}
\newcommand{\new}{\mathrm{new}}
\newcommand{\trial}{{\mathrm{trial}}}

\newcommand{\red}[1]{\textcolor{red}{#1}}
\newcommand{\code}[1]{\texttt{#1}}
\newcommand{\pkg}[1]{\textbf{#1}}
\newcommand{\proglang}[1]{#1}
\newcommand{\citep}[1]{\cite{#1}}

\newcounter{algocf}[section]
\renewcommand{\thealgocf}{\arabic{section}.\arabic{algocf}}
\newcommand{\algo}[4]{
  \refstepcounter{algocf}
  \noindent
  \rule{\textwidth}{0.025cm} \\
  {\bf Algorithm \thealgocf: #2} \\[-0.2cm]
  \rule{\textwidth}{0.025cm} \\[-2\topsep]
  \label{#1}
  \begin{description}
  \item[Input.] #3
  \end{description}
  \begin{enumerate}
    \renewcommand{\labelenumi}{\textbf{\theenumi}.}
    \renewcommand{\labelenumii}{\textbf{\theenumii}.}
    \renewcommand{\labelenumiii}{\theenumiii.}
    \renewcommand{\labelenumiv}{\theenumiv.}
  
    \renewcommand{\theenumi}{Step~\arabic{enumi}}
    \renewcommand{\theenumii}{Step~\arabic{enumi}.\arabic{enumii}}
    \renewcommand{\theenumiii}{\arabic{enumi}.\arabic{enumii}.\arabic{enumiii}}
    \renewcommand{\theenumiv}{\arabic{enumi}.\arabic{enumii}.\arabic{enumiii}.\arabic{enumiv}}
  
    \setlength{\itemsep}{2ex}
    \setlength{\parskip}{0pt}
    \setlength{\parsep}{0pt}
    #4
  \end{enumerate}
  \vspace{-0.5cm}
  \rule{\textwidth}{0.025cm}
}

\begin{document}

\title{Accelerated derivative-free spectral residual method for
  nonlinear systems of equations\thanks{This work was supported by
    FAPESP (grants 2013/07375-0, 2016/01860-1, and 2018/24293-0) and
    CNPq (grants 302538/2019-4 and 302682/2019-8).}}

\author{
  E. G. Birgin\thanks{Department of Computer Science, Institute of
    Mathematics and Statistics, University of S\~ao Paulo, Rua do
    Mat\~ao, 1010, Cidade Universit\'aria, 05508-090, S\~ao Paulo, SP,
    Brazil. e-mail: egbirgin@ime.usp.br, diaulas@ime.usp.br.}
  \and
  J. L. Gardenghi\thanks{Faculty UnB Gama, University of Bras\'ilia,
    \'Area Especial de Ind\'ustria Proje\c{c}\~ao A, Setor Leste,
    Gama, 72444-240, Bras\'ilia, DF. e-mail: john.gardenghi@unb.br.}
  \and
  D. S. Marcondes\footnotemark[2]
  \and
  J. M. Mart\'{\i}nez\thanks{Department of Applied Mathematics,
    Institute of Mathematics, Statistics, and Scientific Computing
    (IMECC), State University of Campinas, 13083-859 Campinas SP,
    Brazil. e-mail: martinez@ime.unicamp.br.}}

\date{April 26, 2021}

\maketitle

\begin{abstract}
Spectral residual methods are powerful tools for solving nonlinear
systems of equations without derivatives. In a recent paper, it was
shown that an acceleration technique based on the Sequential Secant
Method can greatly improve its efficiency and robustness. In the
present work, an \proglang{R} implementation of the method is
presented. Numerical experiments with a widely used test bed compares
the presented approach with its plain (i.e.\ non-accelerated) version
that makes part of the \proglang{R} package \pkg{BB}. Additional
numerical experiments compare the proposed method with \pkg{NITSOL}, a
state-of-the-art solver for nonlinear systems. The comparison shows
that the acceleration process greatly improves the robustness of its
counterpart included in the existent \proglang{R} package.  As a
by-product, an interface is provided between~\proglang{R} and the
consolidated \pkg{CUTEst} collection, which contains over a thousand
nonlinear programming problems of all types and represents a standard
for evaluating the performance of optimization methods. \\


\end{abstract}

\section{Introduction} \label{intro}

Solving nonlinear systems of equations is an ubiquitous problem that appears in a wide range of applied fields such as Physics, Chemistry, Engineering, and Statistics, just to name a few. Moreover, many times, equations are computed using black-box codes and derivatives are not available. Thus, derivative-free solution methods are in order.

Given $F:\R^n \to \R^n$, we consider the problem of finding $x \in \R^n$ such that
\begin{equation} \label{theproblem}
F(x) = 0,
\end{equation}
without making use of derivatives. Observing that~\eqref{theproblem} is equivalent to $x = x - \sigma F(x)$, for any $\sigma>0$, Sequential Residual Methods (namely SANE and DF-SANE) based on the iteration $x^{k+1} = x^k - \sigma_k F(x^k)$, where
\[
\sigma_k = \frac{\|s^{k-1}\|^2}{(y^{k-1})^T s^{k-1}}, \;\;\;\; 
s^{k-1} = x^k - x^{k-1}, \;\; \mbox{ and } \;\; y^{k-1} = F(x^k) - F(x^{k-1}),
\]
were introduced in~\cite{lacruzraydan} and~\cite{lmr}. These methods were inspired by the Barzilai-Borwein step of minimization methods; see \cite{barzilaiborwein,raydan1,raydan2}. Although very popular, in part due to its simplicity, these methods may suffer from slow convergence. On the other hand, their simple and fast iterations made them an adequate choice to provide a global convergent framework to the Sequential Secant approach~\citep{barnes,wolfe}. This choice was explored in~\cite{bmdiis}, where the Accelerated DF-SANE method was introduced. Numerical experiments in~\cite{bmdiis} shown that Accelerated DF-SANE compares favorably to the classical truncated Newton approach implemented in the package \pkg{NITSOL}~\citep{nitsol}, when applied to large-scale problems coming from the discretization of partial differential equations.

In the present work, an \proglang{R}~\citep{r} implementation of Accelerated DF-SANE is introduced. Numerical experiments in~\cite{bmdiis} are complemented with numerical experiments using the widely-used testing environment for optimization \pkg{CUTEst}~\citep{cutest}. Problems in the \pkg{CUTEst} collection are given in SIF (Standard Input Format; see~\cite[Chapters 2 and 7]{lancelot}) and a decoder named SifDec translates the problem into \proglang{Fortran} routines. Therefore, in order to be able to use the \pkg{CUTEst} collection, an interface with the \proglang{R} language is required. Such interface is introduced in the present work; and the authors hope that 
its dissemination in the \proglang{R} community could help in testing and assessing the performance of optimization methods developed in \proglang{R}. Classical sets of problems, like the ones introduced in~\cite{mgh} and~\cite{hock1,hock2}, are included in the \pkg{CUTEst} collection. In addition to the comparison with \pkg{NITSOL}, a comparison with the DF-SANE method implemented within the \pkg{BB} package~\citep{ravi} implemented in \proglang{R} is also provided. 

The rest of this work is organized as follows. The Accelerated DF-SANE method and its convergence theoretical results are condensed in Section~\ref{algandconv}. The \proglang{R} implementation of the method and its usage are described in Section~\ref{implementation}. Numerical results are reported in Section~\ref{experiments}. Conclusions are given in the last section.

\section{Accelerated DF-SANE} \label{algandconv}

In this section, the Accelerated DF-SANE method introduced in~\cite{bmdiis} and its theoretical convergence results are summarized. Roughly speaking, Accelerated DF-SANE performs a nonmonotone line search along the direction of the residue. As a result of a double backtracking, at each iteration~$k$, a trial point $x^{k+1}_{\trial}$ is first computed. Before deciding whether this trial point will be the next iterate~$x^{k+1}$ or not (as it would be the case in the plain DF-SANE in which acceleration is not performed), an accelerated point~$x^{k+1}_{\acc}$ is computed. Following sequential secant ideas, $x^{k+1}_{\acc}$ is given by $x^{k+1}_{\acc} = x^k - S_k  Y_k^\dagger F(x^k)$, where $p>1$ is a given parameter, $\kbar = \max\{ 0, k - p + 1 \}$,
\[
\begin{array}{rcl}
s^j &=& x^{j+1} - x^j \mbox{ for } j = \kbar, \dots, k-1,\\[2mm]
y^j &=& F(x^{j+1}) - F(x^j) \mbox{ for } j = \kbar, \dots, k-1,\\[2mm]
s^k &=& x^{k+1}_{\trial} - x^k,\\[2mm]
y^k &=& F(x^{k+1}_{\trial}) - F(x^k),\\[2mm]
S_k &=& (s^{\kbar}, \dots, s^{k-1}, s^k),\\[2mm]
Y_k &=& (y^{\kbar}, \dots, y^{k-1}, y^k),
\end{array}
\]
and $Y_k^\dagger$ is the Moore-Penrose pseudoinverse of $Y_k$. Then, if $\|F(x^{k+1}_{\acc})\|_2^2 < \|F(x^{k+1}_{\trial})\|_2^2$, the methods defines $x^{k+1} = x^{k+1}_{\acc}$; while $x^{k+1} = x^{k+1}_{\trial}$ in the other case. In practice, $x^{k+1}_{\acc}$ is computed by first finding the minimum norm least-squares solution~$\bar \nu$ of the linear system $Y_k \nu = F(x^{k+1}_{\trial})$ and then defining $x^{k+1}_{\acc} = x^{x+1}_{\trial} - S_k \bar \nu$. The minimum-norm least-squares solution~$\bar \nu$ is computed with a complete orthogonalization of~$Y_k$. The key point is that matrix~$Y_k$ corresponds to remove one column and add one column to matrix~$Y_{k-1}$, keeping the cost of each iteration low; see~\cite[\S5.4]{bmdiis} for details. The whole Accelerated DF-SANE method is given in the algorithm that follows.

\algo{algorithm}{Accelerated DF-SANE}{Let $\gamma \in (0,1)$, $0 < \sigmin < \sigmax < \infty$, $0 < \taumin < \taumax < 1$, positive integers~$M$ and~$p$, a sequence $\{\eta_k\}$ such that $\eta_k > 0$ for all $k \in \N$ and $\lim_{k \to \infty} \eta_k = 0$, and $x_0 \in \R^n$ be given. Set $k \leftarrow 0$.}{

\item \label{algstop} If $F(x^k) = 0$, then terminate the execution of the algorithm.

\item \label{algbl} Choose $\sigma_k$ such that $|\sigma_k| \in [\sigma_{\min},
  \sigma_{\max}]$ and $v^k \in \R^n$ such that $\|v^k\| =
  \|F(x^k)\|$. Compute
  \begin{equation} \label{barfk}
    \bar f_k = \max\{f(x^k), \dots, f(x^{\max\{0,k-M+1\}}) \}.
  \end{equation}
  \begin{enumerate}
    \item Set $\alpha_+ \leftarrow 1$ and $\alpha_- \leftarrow 1$.
    \item Set $d \leftarrow -\sigma_k v^k$ and $\alpha
      \leftarrow \alpha_+$. Consider
      \begin{equation} \label{algarmijo}
        f(x^k + \alpha d) \leq \bar f_k + \eta_k - \gamma \alpha^2 f(x^k).
      \end{equation}
      If (\ref{algarmijo}) holds, then define $d^k = d$ and $\alpha_k
      = \alpha$ and go to \ref{algxtrial}.
    \item Set $d \leftarrow \sigma_k v^k$ and $\alpha
      \leftarrow \alpha_-$. If (\ref{algarmijo}) holds, then define
      $d^k = d$ and $\alpha_k = \alpha$ and go to \ref{algxtrial}.
    \item Choose $\alpha_+^{\new} \in [\taumin \alpha_+,
      \taumax \alpha_+]$ and $\alpha_-^{\new} \in [\taumin \alpha_-,
      \taumax \alpha_{-}]$, set $\alpha_+ \leftarrow \alpha_+^{\new}$,
      $\alpha_- \leftarrow \alpha_-^{\new}$, and go to \ref{algbl}.2.
  \end{enumerate}
  
\item \label{algxtrial} Define $x^{k+1}_{\trial} = x^k + \alpha_k d^k$.

\item \label{algacc} Define $x^{k+1}_{\acc} = x^k - S_k  Y_k^\dagger F(x^k)$, where $\kbar = \max\{ 0, k - p + 1 \}$,
    \[
    \begin{array}{rcl}
      s^j &=& x^{j+1} - x^j \mbox{ for } j = \kbar, \dots, k-1,\\[2mm]
      y^j &=& F(x^{j+1}) - F(x^j) \mbox{ for } j = \kbar, \dots, k-1,\\[2mm]
      s^k &=& x^{k+1}_{\trial} - x^k,\\[2mm]
      y^k &=& F(x^{k+1}_{\trial}) - F(x^k),\\[2mm]
      S_k &=& (s^{\kbar}, \dots, s^{k-1}, s^k),\\[2mm]
      Y_k &=& (y^{\kbar}, \dots, y^{k-1}, y^k),
    \end{array}
    \]
    and $Y_k^\dagger$ is the Moore-Penrose pseudoinverse of $Y_k$.

\item Choose $x^{k+1} \in \left\{ x^{k+1}_{\trial}, x^{k+1}_{\acc}
    \right\}$ such that
    \[
    \|F(x^{k+1})\| = \min \left\{ \|F(x^{k+1}_{\trial})\|,
    \|F(x^{k+1}_{\acc})\| \right\}.
    \]

\item Set $k \leftarrow k+1$, and go to \ref{algstop}.
}

In practice, at \ref{algstop}, given $\varepsilon > 0$, the stopping criterion $\|F(x^k)\|=0$ is replaced with
\begin{equation} \label{scsuccess}
\|F(x^k)\|_2 \leq \varepsilon.
\end{equation} 
(Criterion $\|F(x^k)\|=0$ in the algorithm is necessary so we can state theoretical asymptotic properties of an infinite sequence generated by the algorithm.) At \ref{algbl}, the spectral choice for~$\sigma_k$ (see~\cite{barzilaiborwein,raydan1,raydan2,bmr,bmracmtoms,bmrency,bmr5}) corresponds to
\[
\sigma_k^{\spg} = \frac{(x^k - x^{k-1})^T(x^k -
  x^{k-1})}{(x^k - x^{k-1})^T(F(x^k) - F(x^{k-1})}.
\]
Following~\cite{lmr}, if $|\sigma_k^{\spg}| \in [\sigma_{\min},\min\{1,\sigma_{\max}\}]$, then we take $\sigma_k = \sigma_k^{\spg}$; otherwise, we take $\sigma_k = \max\{ \sigma_{\min}, \min\{ \| x^k \|_2 / \| v^k \|_2, \sigma_{\max} \} \}$. Still at \ref{algbl}, the residual choice for the search direction corresponds to $v_k=F(x^k)$. At \ref{algbl}.4, we compute $\alpha_+^{\new}$ as the minimizer of the univariate quadratic $q(\alpha)$ that interpolates $q(0)=f(x^k)$, $q(\alpha_+)=f(x^k-\alpha_+ \sigma_k F(x^k))$, and $q'(0)=-\sigma_k F(x^k)^T \nabla f(x^k) = -\sigma_k F(x^k)^T J(x^k) F(x^k)$. Following~\cite{lmr}, since we consider $J(x^k)$ unavailable, we consider $J(x^k) = I$. Thus,
\[
\alpha_+^{\new} = \max \left\{ \taumin \alpha_+, \min \left\{
\frac{\alpha_+^2 f(x^k)}{f(x^k - \alpha_+ \sigma_k F(x^k)) +
  (2\alpha_+-1) f(x^k)}, \taumax \alpha_+ \right\} \right\}.
\]
Analogously,
\[
\alpha_-^{\new} = \max \left\{ \taumin \alpha_-, \min \left\{
\frac{\alpha_-^2 f(x^k)}{f(x^k + \alpha_- \sigma_k F(x^k)) +
  (2\alpha_--1) f(x^k)}, \taumax \alpha_- \right\} \right\}.
\]

Theoretical results of Algorithm~\ref{algorithm} are given in~\cite[\S3 and \S4]{bmdiis}. Briefly, limit points of sequences generated by the algorithm are solutions of the nonlinear system or the gradient of the corresponding sum of squares is null. Moreover, under suitable assumptions, the convergence to solutions is superlinear.

\section{Usage of the R implementation} \label{implementation}

We implemented Algorithm~\ref{algorithm} in \proglang{R} language as a subroutine named \code{dfsaneacc}. Codes are freely available at \url{https://github.com/johngardenghi/dfsaneacc} and at the Journal web page accompanying the present work. In this section, we describe how to use \code{dfsaneacc} to solve a nonlinear system coded in \proglang{R} and how to solve a nonlinear system from the \pkg{CUTEst} collection.

The calling sequence of \code{dfsaneacc} is given by
\begin{verbatim}
  R> dfsaneacc(x, evalr, nhlim, epsf, maxit, iprint, ...)
\end{verbatim}
where
\begin{description}
\item[\code{x}:] is an $n$-dimensional array containing the initial guess.
\item[\code{evalr}:] is the subroutine that computes $F$ at a point \code{x}. This subroutine must have the calling sequence
\begin{verbatim}
  evalr <- function(x, ...) {}
\end{verbatim}
where \code{...} represents the additional arguments of \code{dfsaneacc}. The subroutine must return $F$ evaluated at \code{x}.
\item[\code{nhlim}:] corresponds to $p+1$, where $p \geq 1$ is the integer that says how many previous iterates must be considered in the Sequential Secant acceleration at \ref{algacc}. The ``default'' value is $p=5$, so \code{nhlim}=6; but having a problem at hand, it is recommendable to try different values.
\item[\code{epsf}:] corresponds to the stopping tolerance $\varepsilon$ in~\eqref{scsuccess}.
\item[\code{maxit}:] represents the maximum number of iterations. It default value is \code{maxit}=$+\infty$.
\item[\code{iprint}:] determines the level of the details in the output of the routine -- \code{iprint}=$-1$ means no output, \code{iprint}=$0$ means basic information at every iteration, \code{iprint}=$1$ adds additional information related to the backtracking strategy (\ref{algbl}), and \code{iprint}=$2$ adds information related to the computation of the acceleration step (\ref{algacc}). Its default value is \code{iprint}=$-1$.
\end{description}

As an example, consider the \textit{Exponential Function 2} from~\cite[p.596]{lacruzraydan} given by $F(x) = \left( F_1(x), \dots, F_n(x) \right)^T$, where
\setlength{\arraycolsep}{2pt}
\[
\begin{array}{rcl}
F_1(x) & = & \displaystyle e^{x_1} - 1 \\[1ex]
F_i(x) & = & \displaystyle \frac{i}{10} (e^{x_1} + x_{i-1} - 1) \mbox{ for } i = 2, \dots, n,
\end{array}
\]
with the initial guess $x^0 = (\frac{1}{n^2},\dots,\frac{1}{n^2})^T$. The first step is to code it in \proglang{R} as follows:
\begin{verbatim}
  R> expfun2 <- function(x) {
  +    n <- length(x)
  +    f <- rep(NA, n)
  +    f[1] <- exp(x[1]) - 1.0
  +    f[2:n] <- (2:n)/10.0 * (exp(x[2:n]) + x[1:n-1] - 1)
  +    f
  +  }
\end{verbatim}
Then, we set the dimension $n$ and the initial point~$x^0$ and call \code{dfsaneacc} as follows:
\begin{verbatim}
  R> n <- 3
  R> x0 <- rep(1/n^2, n)
  R> ret <- dfsaneacc(x=x0, evalr=expfun2, nhlim=6, epsf=1.0e-6*sqrt(n), 
  +                   iprint=0)
\end{verbatim}
obtaining the result below:
\begin{verbatim}
  Iter:  0  f =  0.02060606
  Iter:  1  f =  0.001215612
  Iter:  2  f =  4.68925e-05
  Iter:  3  f =  4.654419e-08
  Iter:  4  f =  1.135198e-11
  Iter:  5  f =  9.154603e-16
  success!

  $x
                [,1]
  [1,] -3.582692e-11
  [2,] -7.222425e-08
  [3,] -1.638214e-08

  $res
  [1] -3.582690e-11 -1.445201e-08 -2.658192e-08

  $normF
  [1] 9.154603e-16

  $iter
  [1] 5

  $fcnt
  [1] 11

  $istop
  [1] 0
\end{verbatim}
where
\begin{description}
\item[\code{x}:] is the approximation to a solution $x_*$.
\item[\code{res}:] corresponds to $F(x_*)$.
\item[\code{normF}:] corresponds to $f(x_*)=\|F(x_*)\|_2^2$.
\item[\code{iter}:] is the number of iterations.
\item[\code{fcnt}:] is the number of calls to \code{evalr}, i.e.\ the number of functional evaluations.
\item[\code{istop}:] is the exit code, where \code{istop}=0 means that~$x_*$ satisfies~\eqref{scsuccess}, i.e.\ $\|F(x_*)\|_2 \leq \varepsilon$, and \code{istop}=1 means that the maximum allowed number of iterations was reached.
\end{description}

In the rest of this section, we show how to solve a nonlinear system from the \pkg{CUTEst} collection. \pkg{CUTEst} can be downloaded from \url{https://github.com/ralna/CUTEst}. It is assumed that \pkg{CUTEst} is installed, in particular \code{SifDec}, and that there is a folder with all problems in SIF format. 

The first step is to choose a problem and run \code{SifDec} that, based on the problem's SIF file, generates a \proglang{Fortran} routine to evaluate, in this case, function~$F$. It should be mentioned that problems in the \pkg{CUTEst} collection are general nonlinear optimization problems of the form
\begin{equation} \label{cutestprob}
\mbox{Minimize } \Phi(x) \mbox{ subject to } h(x)=0, \; \ell_g \leq g(x) \leq u_g, \; \ell \leq x \leq u,
\end{equation}
where $\Phi : \R^n \to \R$ is the objective function, $h:\R^n \to \R^{m_E}$ represents $m_E$ equality constraints, $g:\R^n \to \R^{m_I}$ represents $m_I$ two-side inequality constraints, $\ell_g, u_g \in \R^{m_I}$, and $\ell, u \in \R^n$ represent bounds on the variables. (Some components of $\ell_g$ and $\ell$ can be $-\infty$ as well as some components of $u_g$ and $u$ can be equal to $+\infty$.) Thus, a nonlinear system of equations corresponds to a problem of the form~\eqref{cutestprob} with constant or null objective function, equality constraints only, and $n=m_E$; and, in the context of the present work, we define $F(x) \equiv h(x)$. Once the \proglang{Fortran} codes have been generated, a dynamic library must be built and loaded in~\proglang{R}. The wrapper (written in \proglang{R}) uses this library to call, using the \code{.Call} tool, a \proglang{C} subroutine from an existent \proglang{C} interface of \pkg{CUTEst}, that calls the generated \proglang{Fortran} subroutine. In fact, \pkg{CUTEst} is mainly coded in \proglang{Fortran} and calling a \proglang{Fortran} subroutine using the tool \code{.Fortran} would be the natural choice. However, numerical experiments shown that the combination of \code{.Call} with the existent \proglang{C} interface of \pkg{CUTEst} is faster.

The wrapper consists in five routines named \code{cutest\_init}, \code{cutest\_end}, \code{cutest\_getn}, \code{cutest\_getx0}, and \code{cutest\_evalr}. Routine \code{cutest\_init} receives as parameter the name of a problem and executes all initialization tasks described in the previous paragraph. Routine \code{cutest\_end} has no parameters and it cleans the environment by freeing the memory allocated in the call to \code{cutest\_init}. The other three routines are self-explanatory. So, for example, a problem named \textsc{Booth} can be solved simply by typing:
\begin{verbatim}
  R> cutest_init('BOOTH')
  R> n <- cutest_getn()
  R> x0 <- cutest_getx0()
  R> ret <- dfsaneacc(x=x0, evalr=cutest_evalr, nhlim=6, epsf=1.0e-6*sqrt(n),
     +                iprint=0)
  R> cutest_end()
\end{verbatim}
The output follows:
\begin{verbatim}
  Iter:  0  f =  74
  Iter:  1  f =  3.544615
  Iter:  2  f =  9.860761e-31
  success!

  $x
       [,1]
  [1,]    1
  [2,]    3

  $res
  [1] -8.881784e-16 -4.440892e-16

  $normF
  [1] 9.860761e-31

  $iter
  [1] 2

  $fcnt
  [1] 7

  $istop
  [1] 0
\end{verbatim}

There are environment variables that must be set to indicate where \pkg{CUTEst} was installed, which is the folder that contains the SIF files of the problems, and which \proglang{Fortran} compiler and compiling options must be used. A README file with detailed instructions accompanies the distribution of Accelerated DF-SANE and the \pkg{CUTEst} interface with~\proglang{R}.

\section{Numerical experiments} \label{experiments}

In this section, we show the performance of Algorithm~\ref{algorithm} by putting it in perspective in relation to the DF-SANE algorithm of the \pkg{BB} package~\citep{ravi} and the well-known
 \pkg{NITSOL} method~\citep{nitsol}. For that, we consider \textit{all}~70 nonlinear systems of the \pkg{CUTEst} collection~\citep{cutest} with their default dimensions and their default initial points.

In this work, we implemented Algorithm~\ref{algorithm} in~\proglang{R}; while a \proglang{Fortran} implementation, available at \url{https://www.ime.usp.br/~egbirgin/sources/accelerated-df-sane/}, was given in~\cite{bmdiis}. The state-of-the-art solver \pkg{NITSOL} is available in \proglang{Fortran} in \url{https://users.wpi.edu/~walker/NITSOL/}. A \proglang{Fortran} version of DF-SANE is available under request to the authors of~\cite{lmr}; while an \proglang{R} implementation of DF-SANE is available as part of the \pkg{BB} package~\citep{ravi}. Problems of the \pkg{CUTEst} collection are written in SIF (Standard Input Format); and a tool named SifDec (SIF Decoder) generates \proglang{Fortran} routines to evaluate the objective function, in addition to constraints and their derivatives when desired. So, an interface between \proglang{R} and \pkg{CUTEst} was implemented in order to test DF-SANE and Accelerated DF-SANE (both in \proglang{R}) with the problems of the \pkg{CUTEst} collection. \proglang{Fortran} codes were compiled with the GFortran compiler of GCC (version 9.3.0). \proglang{R} codes were run in version~4.0.2. Tests were conducted on a computer with an Intel Core i7 7500 processor and 12 GB of RAM memory, running Linux (Ubuntu 20.10). 

Regarding the DF-SANE method~\citep{lmr} that is available as part of the \pkg{BB} package~\citep{ravi}, a few considerations are in order. First of all, in the numerical experiments, we considered function \code{dfsane} from package \pkg{BB} version 2019.10-1. In the \pkg{BB} package, there is a routine named \code{BBsolve} that is a wrapper for \code{dfsane}. \code{BBsolve} calls \code{dfsane} repeatedly with different algorithm parameters aiming to find a solution to the problem at hand. Since this strategy can be used in connection with any method, aiming a fair comparison, in the present work we report the results obtained with a single run of \code{dfsane} with its default parameters. This means that the strategies described in~\cite[\S2.4]{ravi} are not being considered. On the other hand, \code{dfsane} improves the original DF-SANE method introduced in~\cite{lmr} in several ways; see~\cite[\S2.3]{ravi}. Among the improvements, there is one that is particularly relevant in the context of the present work: when the plain DF-SANE method fails by lack of progress, \code{dfsane} launches an alternative method -- it runs L-BFGS-B for the minimization of $f(x)=\|F(x)\|_2^2$. L-BFGS-B~\citep{lbfgsb} is a limited-memory quasi-Newton method for bound-constrained minimization. In some way, it could be said that this modification aims to mitigate the slow convergence of DF-SANE. In contrast to the approach presented in the present paper, this device is triggered only once slow convergence has been detected; while in the present work, acceleration is done at every iteration. Anyway, it is worth noticing that, by comparing the method being introduced in the present work with \code{dfsane} from the \pkg{BB} package, a comparison is being done with an improved version of the original DF-SANE introduced in~\cite{lmr}.

From now on, we refer to the DF-SANE of the \pkg{BB} package simply as DF-SANE; while we refer to Algorithm~\ref{algorithm} as ``Accelerated DF-SANE''. \pkg{NITSOL} includes three main iterative solvers for linear systems: GMRES, BiCGSTAB, and TFQMR. Numerical experiments showed that, on the considered set of problems, using GMRES presents the best performance among the three options. So, from now on, we refer to \pkg{NITSOL} as ``\pkg{NITSOL} (GMRES)''. All default parameters of DF-SANE and \pkg{NITSOL} (GMRES) were considered. For the Accelerated DF-SANE, following~\cite{bmdiis}, we considered $\gamma=10^{-4}$, $\taumin=0.1$, $\taumax=0.5$, $M=10$, $\sigmin=\sqrt{\epsilon}$, $\sigmax=1/\sqrt{\epsilon}$, $\eta_k = 2^{-k} \min\{ \frac{1}{2} \| F(x^0)\|, \sqrt{\| F(x^0) \| }\}$, where $\epsilon \approx 10^{-16}$ is the machine precision, and $p=5$. To promote a fair comparison, in all three methods, the common stopping criterion~\eqref{scsuccess} with $\varepsilon = 10^{-6} \sqrt{n}$, was considered. In addition, each method has its own alternative stopping criteria, mainly related to lack of progress; and a CPU time limit of 3 minutes per method/problem was also imposed in the numerical experiments.

Table~\ref{accvsravi} shows the result of DF-SANE and Accelerated DF-SANE (recall that both methods are coded in~\proglang{R}). In the table, the first two columns show the problem name and the number of variables and equations. Then, for each method, the table reports the value of $\|F(x)\|_2$ at the final iterate (column $\|F(x_*)\|_2$), the number of iterations (column \#iter), the number of functional evaluations (column \#feval), and the CPU time in seconds (column time). In column $\|F(x_*)\|_2$, figures in red are the ones that \textit{do not} satisfy~\eqref{scsuccess}. It is worth noticing that in all cases in which the final iterate of DF-SANE does not satisfy~\eqref{scsuccess}, DF-SANE stops by ``lack of progress'' (flag equal to~5). When the same happens with Accelerated DF-SANE, since no stopping criterion due to lack of progress was implemented, it stops by reaching the CPU time limit. The table shows that Accelerated DF-SANE satisfied the stopping criterion~\eqref{scsuccess} related to success in~44 out of the 70 considered problems; while DF-SANE did the same in~32 problems. Moreover, there were 30 problems that were solved by both methods, 14 problems that were solved by Accelerated DF-SANE only, and 2 problems that were solved by DF-SANE only. These figures show that the acceleration step improves the robustness of DF-SANE.

\LTcapwidth=\textwidth
\begin{center}
{\footnotesize\tabcolsep=3pt
\begin{longtable}{lr@{\hskip 2ex}rrrrrrrrr}
\toprule
\multirow{2}{*}{Problem} &
\multirow{2}{*}{$n$} &
\multicolumn{4}{c}{Accelerated DF-SANE} & &
\multicolumn{4}{c}{DF-SANE}\\ 
\cmidrule{3-6} \cmidrule{8-11}
&        & $\|F(x_*)\|$&  \# iter  &  \# feval   &     time     && $\|F(x_*)\|$&  \# iter  &  \# feval   &     time \\ \midrule
     BOOTH &      2 &       $9.9$E$-16$ &          2 &          7 &   0.005065 &&       $2.4$E$-07$ &          7 &          8 &   0.004709 \\  
   CLUSTER &      2 &       $8.3$E$-07$ &         23 &        108 &   0.007488 &&       $2.4$E$-07$ &         40 &         42 &   0.005575 \\  
    CUBENE &      2 &       $4.0$E$-13$ &          9 &         20 &   0.005451 &&       $1.0$E$-06$ &         24 &         26 &   0.005079 \\  
DENSCHNCNE &      2 &       $2.3$E$-11$ &         10 &         23 &   0.005626 &&       $1.4$E$-07$ &         16 &         17 &   0.005019 \\  
DENSCHNFNE &      2 &       $2.7$E$-07$ &          7 &         23 &   0.005479 &&       $2.5$E$-07$ &         27 &         40 &   0.005340 \\  
  FREURONE &      2 &       $1.5$E$-08$ &         16 &         55 &   0.006213 && \red{$1.1$E$+01$} &        103 &        123 &   0.007488 \\  
    GOTTFR &      2 &       $1.3$E$-07$ &         23 &         67 &   0.006572 && \red{$2.6$E$-02$} &      24196 &     154606 &   2.629654 \\  
  HIMMELBA &      2 &       $0.0$E$+00$ &          2 &          7 &   0.005166 &&       $1.3$E$-07$ &          7 &          8 &   0.004738 \\  
  HIMMELBC &      2 &       $8.4$E$-08$ &          5 &         13 &   0.005269 &&       $7.0$E$-07$ &         10 &         11 &   0.004874 \\  
  HIMMELBD &      2 & \red{$2.4$E$+00$} &     211279 &    5989534 & 180.000000 && \red{$2.4$E$+00$} &        188 &        211 &   0.009555 \\  
       HS8 &      2 &       $4.4$E$-08$ &          5 &         13 &   0.005335 &&       $2.1$E$-07$ &         14 &         15 &   0.004822 \\  
    HYPCIR &      2 &       $8.7$E$-10$ &          6 &         14 &   0.005353 &&       $1.2$E$-06$ &         13 &         14 &   0.004824 \\  
  POWELLBS &      2 & \red{$2.3$E$-03$} &     225728 &    4561326 & 180.000000 &&       $8.4$E$-07$ &        106 &        367 &   0.010667 \\  
  POWELLSQ &      2 & \red{$3.9$E$+00$} &     317171 &     779427 & 180.000000 && \red{$9.8$E$-03$} &     665188 &    6522441 & 101.318056 \\  
  PRICE3NE &      2 &       $3.9$E$-10$ &          7 &         19 &   0.005414 &&       $9.0$E$-07$ &         15 &         16 &   0.004841 \\  
  PRICE4NE &      2 &       $1.3$E$-10$ &         10 &         27 &   0.005625 &&       $2.0$E$-08$ &         37 &         39 &   0.005394 \\  
   RSNBRNE &      2 &       $4.4$E$-16$ &         56 &        204 &   0.009382 &&       $3.7$E$-07$ &        428 &        564 &   0.018345 \\  
  SINVALNE &      2 &       $4.9$E$-15$ &         16 &         77 &   0.006542 && \red{$2.1$E$+00$} &       5063 &      52078 &   0.846892 \\  
 WAYSEA1NE &      2 &       $1.3$E$-10$ &         12 &         36 &   0.005866 &&       $1.0$E$-06$ &        785 &       3466 &   0.065970 \\  
 WAYSEA2NE &      2 &       $8.4$E$-07$ &        481 &       2179 &   0.052801 && \red{$3.4$E$+01$} &     714039 &   12109386 & 180.004208 \\  
DENSCHNDNE &      3 &       $2.1$E$-07$ &         26 &         62 &   0.006747 &&       $1.1$E$-06$ &         83 &         86 &   0.006548 \\  
DENSCHNENE &      3 &       $9.6$E$-11$ &          6 &         16 &   0.005380 && \red{$9.8$E$-01$} &        107 &        112 &   0.007418 \\  
   HATFLDF &      3 &       $1.4$E$-08$ &         26 &         78 &   0.006926 &&       $9.6$E$-07$ &        586 &        907 &   0.024690 \\  
HATFLDFLNE &      3 & \red{$8.0$E$-03$} &     216456 &    5660213 & 180.000000 && \red{$8.2$E$-03$} &        170 &        251 &   0.010048 \\  
   HELIXNE &      3 &       $2.8$E$-09$ &         13 &         35 &   0.005898 && \red{$3.1$E$+01$} &        102 &        574 &   0.013981 \\  
  HIMMELBE &      3 &       $1.2$E$-15$ &          9 &         21 &   0.005566 && \red{$2.1$E$+00$} &        127 &        128 &   0.007795 \\  
    RECIPE &      3 &       $2.9$E$-07$ &         58 &        355 &   0.012444 &&       $1.4$E$-06$ &         56 &         57 &   0.005821 \\  
  ZANGWIL3 &      3 &       $1.4$E$-14$ &          3 &         11 &   0.005174 &&       $1.3$E$-08$ &         25 &         27 &   0.005093 \\  
  POWELLSE &      4 &       $7.3$E$-07$ &         24 &         70 &   0.006980 && \red{$1.5$E$+01$} &        101 &        240 &   0.009092 \\  
POWERSUMNE &      4 & \red{$4.6$E$-03$} &       2761 &      64429 & 180.000000 && \red{$2.0$E$-02$} &        411 &        485 &   0.017388 \\  
    HEART6 &      6 &       $7.2$E$-07$ &     245873 &    3845345 &  67.912296 && \red{$1.9$E$+01$} &        116 &        476 &   0.013026 \\  
    HEART8 &      8 &       $2.2$E$-06$ &      54602 &     823267 &  14.860346 && \red{$1.3$E$+01$} &        101 &        332 &   0.010646 \\  
  COOLHANS &      9 &       $1.5$E$-06$ &         10 &         45 &   0.006065 && \red{$3.5$E$-02$} &        120 &        124 &   0.007696 \\  
  MOREBVNE &     10 &       $1.6$E$-06$ &         37 &        219 &   0.009777 &&       $3.0$E$-06$ &         73 &         76 &   0.006361 \\  
  OSCIPANE &     10 & \red{$1.0$E$+00$} &         54 &        707 & 180.000000 && \red{$1.0$E$+00$} &        100 &        113 &   0.007410 \\  
 TRIGON1NE &     10 &       $1.9$E$-06$ &         13 &         29 &   0.005877 &&       $1.7$E$-06$ &         30 &         33 &   0.005321 \\  
   INTEQNE &     12 &       $9.2$E$-07$ &          3 &          7 &   0.005143 &&       $1.2$E$-06$ &          5 &          6 &   0.004616 \\  
   HATFLDG &     25 &       $5.0$E$-06$ &      13389 &     211286 &   4.406962 && \red{$5.0$E$+00$} &        102 &        189 &   0.008855 \\  
   HYDCAR6 &     29 & \red{$2.3$E$-02$} &     206865 &    4255024 & 180.000000 && \red{$2.5$E$+01$} &        102 &        430 &   0.014045 \\  
  METHANB8 &     31 & \red{$3.9$E$-03$} &     198664 &    4495087 & 180.000000 && \red{$9.9$E$-01$} &        102 &        109 &   0.007866 \\  
  METHANL8 &     31 & \red{$1.6$E$-01$} &     173606 &    3542099 & 180.000000 && \red{$6.5$E$+01$} &        101 &        490 &   0.015252 \\  
  HYDCAR20 &     99 & \red{$2.3$E$-01$} &     170393 &    3142121 & 180.000000 && \red{$3.6$E$+01$} &        101 &        335 &   0.016278 \\  
  LUKSAN21 &    100 &       $8.9$E$-06$ &         48 &        441 &   0.016229 &&       $6.7$E$-06$ &         69 &         88 &   0.006922 \\  
 MANCINONE &    100 &       $5.9$E$-07$ &          5 &         17 &   0.022032 &&       $5.2$E$-06$ &          7 &          8 &   0.012426 \\  
    QINGNE &    100 &       $4.8$E$-06$ &         21 &         45 &   0.006954 &&       $4.5$E$-06$ &         30 &         36 &   0.005532 \\  
   ARGTRIG &    200 &       $1.2$E$-05$ &         57 &        199 &   0.030535 &&       $1.2$E$-05$ &         80 &         87 &   0.014297 \\  
  BROWNALE &    200 &       $1.0$E$-05$ &          9 &         25 &   0.007390 &&       $1.2$E$-07$ &         15 &         16 &   0.005847 \\  
  CHANDHEU &    500 &       $1.4$E$-05$ &         18 &         99 &   0.273017 &&       $2.2$E$-05$ &         95 &        104 &   0.286036 \\  
  10FOLDTR &   1000 & \red{$9.3$E$+06$} &       8222 &     245098 & 180.000000 && \red{$1.8$E$+05$} &        183 &       1167 &   0.845994 \\  
       KSS &   1000 &       $9.3$E$-06$ &          5 &         17 &   0.028989 &&       $7.5$E$-06$ &          9 &         12 &   0.021450 \\  
    MSQRTA &   1024 & \red{$6.1$E$+01$} &      24241 &     454743 & 180.000000 && \red{$8.6$E$+01$} &        129 &        585 &   0.227472 \\  
    MSQRTB &   1024 & \red{$5.7$E$+01$} &      26216 &     450488 & 180.000000 && \red{$8.6$E$+01$} &        123 &        615 &   0.239714 \\  
   EIGENAU &   2550 & \red{$1.7$E$+02$} &       5138 &     103264 & 180.000000 && \red{$1.8$E$+02$} &        118 &        563 &   0.987960 \\  
    EIGENB &   2550 & \red{$9.8$E$+00$} &       6918 &     102189 & 180.000000 && \red{$9.9$E$+00$} &        856 &       7459 &  12.716400 \\  
    EIGENC &   2652 & \red{$1.0$E$+02$} &       4916 &      97087 & 180.000000 && \red{$1.0$E$+02$} &        112 &        545 &   1.014087 \\  
NONMSQRTNE &   4900 & \red{$2.4$E$+02$} &       3252 &      43571 & 180.000000 && \red{$2.2$E$+02$} &       7353 &      47727 & 180.023645 \\  
  BROYDN3D &   5000 &       $5.3$E$-05$ &         12 &         25 &   0.025578 &&       $1.7$E$-05$ &         16 &         17 &   0.010604 \\  
  BROYDNBD &   5000 & \red{$1.0$E$+00$} &      31283 &     472515 & 180.000000 && \red{$3.6$E$+01$} &        124 &        327 &   0.132678 \\  
  BRYBNDNE &   5000 & \red{$1.0$E$+00$} &      31192 &     471278 & 180.000000 && \red{$3.6$E$+01$} &        124 &        327 &   0.132835 \\  
  NONDIANE &   5000 & \red{$1.4$E$+00$} &      33386 &     716126 & 180.000000 && \red{$6.4$E$+02$} &        102 &        483 &   0.129502 \\  
 SBRYBNDNE &   5000 & \red{$2.7$E$+02$} &      18630 &     377758 & 180.000000 && \red{$2.6$E$+02$} &        319 &        897 &   0.356915 \\  
SROSENBRNE &   5000 &       $3.1$E$-09$ &          9 &         34 &   0.020881 &&       $5.7$E$-08$ &         23 &         25 &   0.012307 \\  
SSBRYBNDNE &   5000 & \red{$1.8$E$+02$} &      23551 &     354751 & 180.000000 && \red{$1.3$E$+02$} &        302 &       1192 &   0.460639 \\  
TQUARTICNE &   5000 & \red{$8.7$E$-01$} &      53163 &     550903 & 180.000000 && \red{$8.9$E$-01$} &        790 &       3991 &   0.853161 \\  
  OSCIGRNE & 100000 &       $1.8$E$-04$ &         28 &         66 &   1.013625 &&       $2.0$E$-04$ &         24 &         25 &   0.196684 \\  
   CYCLIC3 & 100002 & \red{$6.8$E$-01$} &       1921 &      27552 & 180.000000 &&       $2.3$E$-04$ &      11410 &      11765 &  83.093461 \\  
  YATP1CNE & 123200 &       $2.6$E$-07$ &         14 &         41 &   1.443373 && \red{$8.4$E$+03$} &        103 &        865 &  20.785781 \\  
   YATP1NE & 123200 &       $2.6$E$-07$ &         14 &         41 &   1.445582 && \red{$8.4$E$+03$} &        103 &        865 &  20.736302 \\  
  YATP2CNE & 123200 & \red{$3.1$E$+04$} &        606 &       8821 & 180.000000 && \red{$7.2$E$+04$} &        114 &        830 &  16.063343 \\  
   YATP2SQ & 123200 & \red{$4.3$E$+04$} &        723 &       8917 & 180.000000 && \red{$4.5$E$+04$} &        104 &        115 &   2.406395 \\  
\bottomrule 
\caption{Detailed results of the application of Accelerated DF-SANE and DF-SANE to the~70 considered problems from the \pkg{CUTEst} collection.}
\label{accvsravi}
\end{longtable}
}
\end{center}

\vspace{-3em}

\begin{figure}[ht!]
\includegraphics[width=0.5\textwidth]{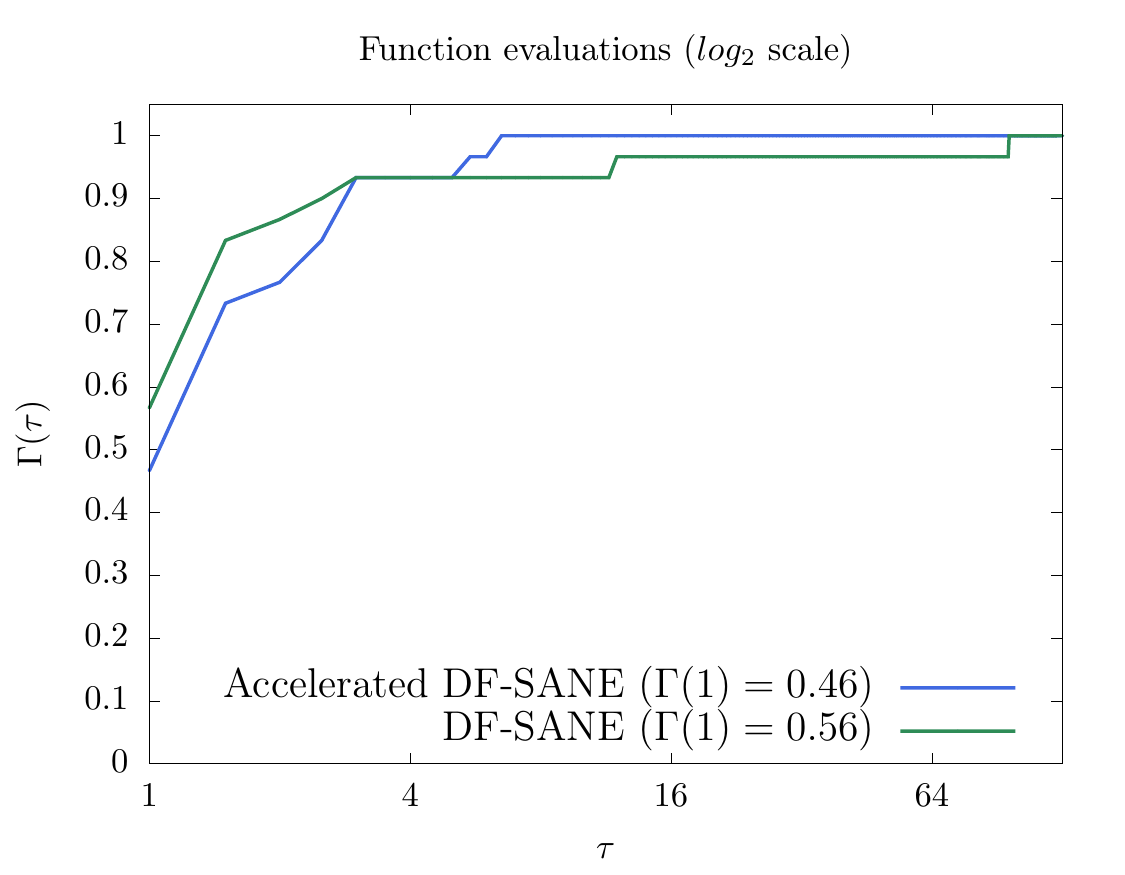}
\includegraphics[width=0.5\textwidth]{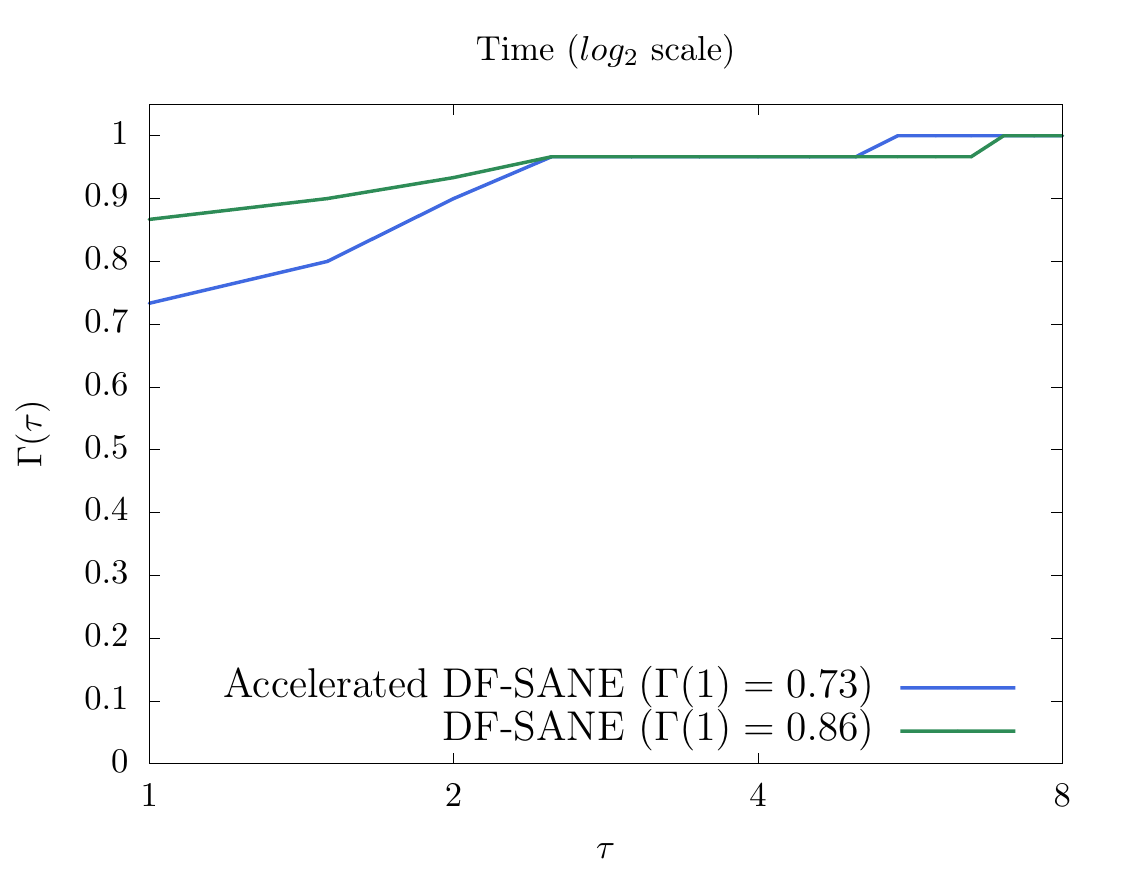}
\caption{Performance profiles of Accelerated DF-SANE and DF-SANE considering the~30 problems from the \pkg{CUTEst} collection in which both methods found a solution.}
\label{ppdiisravi}
\end{figure}

Figure~\ref{ppdiisravi} compares the methods' efficiencies using performance profiles~\citep{pp}. In a performance profile, for $i \in M = \{ \mbox{Accelerated DF-SANE}, \mbox{DF-SANE} \}$,
\[
\Gamma_i(\tau) = \frac{\# \left\{ j \in \{ 1,\dots,n_P\} \;|\; t_{ij} \leq \tau \min_{m \in M} \{ t_{mj} \}\right\}}{n_P},
\]
where $\#{\cal S}$ denotes the cardinality of set ${\cal S}$, $n_P=70$ is the number of problems being considered, and $t_{ij}$ is a measure of the performance of method~$i$ when applied to problem~$j$. If method~$i$ was not able to solve problem~$j$, then we set $t_{ij}=+\infty$. With these definitions, $\Gamma_i(1)$ is the fraction of problems in which method~$i$ was the fastest method to find a solution; while $\Gamma_i(\tau)$ for $\tau$ sufficiently large is the fraction of problems that method~$i$ was able to solve, independently of the required effort. Another possibility, once the robustness of the methods being compared has been established, is to restrict the set of problems in a performance profile to the set of problems that were solved by both methods ($n_P=30$ in this case); so $t_{ij} < +\infty$ for all~$i$ and $j$. With these definitions, the performance profile does not reflect the robustness of the methods any more ($\Gamma_i(\tau)=1$ for a sufficiently large $\tau$ for all $i \in M$) and it is focused on the methods' efficiency. ($\Gamma_i(1)$ still represents the fraction of problems in which method~$i$ was the fastest method to find a solution.) This was the choice in Figure~\ref{ppdiisravi}, in which the number of functional evaluations and the CPU time were used as performance measures. Both graphics show the methods have very similar efficiencies. It is worth noticing that CPU times smaller than 0.01 seconds are considered as being 0.01 and that approximately 90\% of the CPU times, associated with the problems that both methods solve, are smaller than 0.1 seconds.

\begin{figure}[h]
\includegraphics[width=0.5\textwidth]{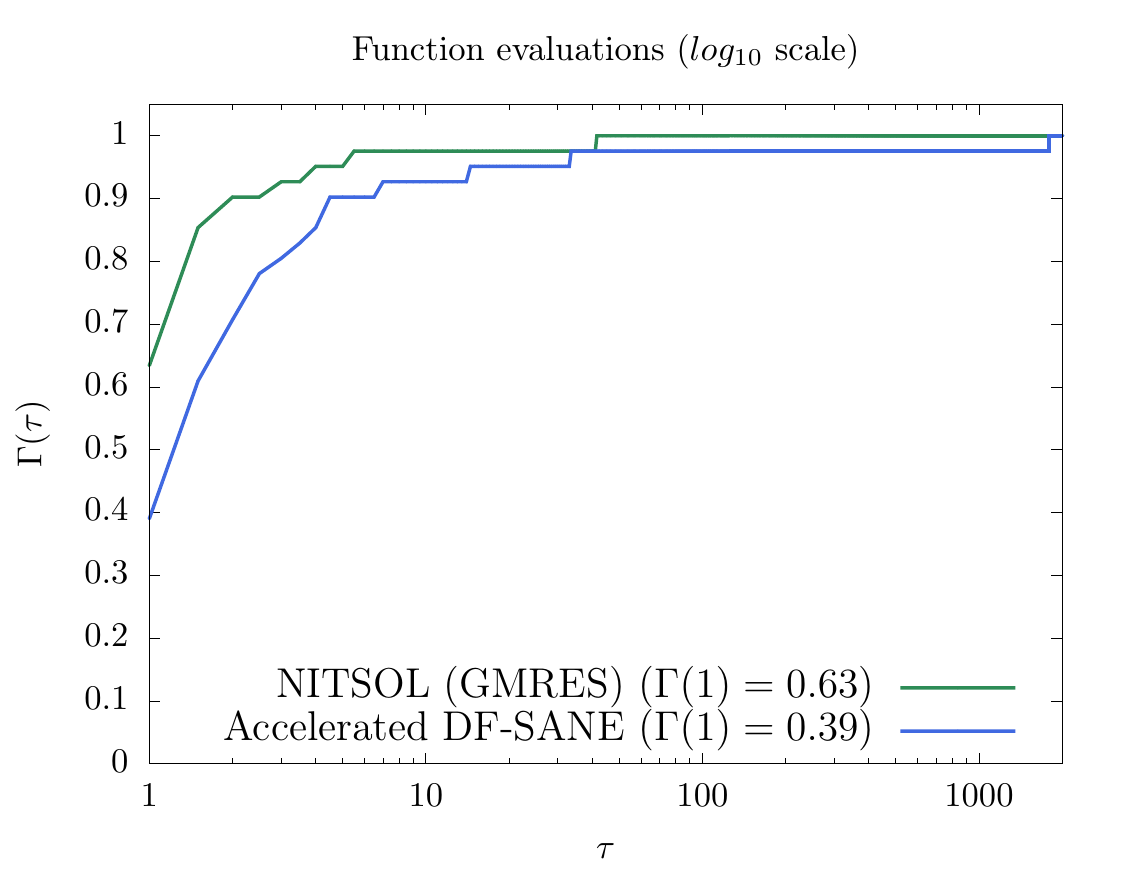}
\includegraphics[width=0.5\textwidth]{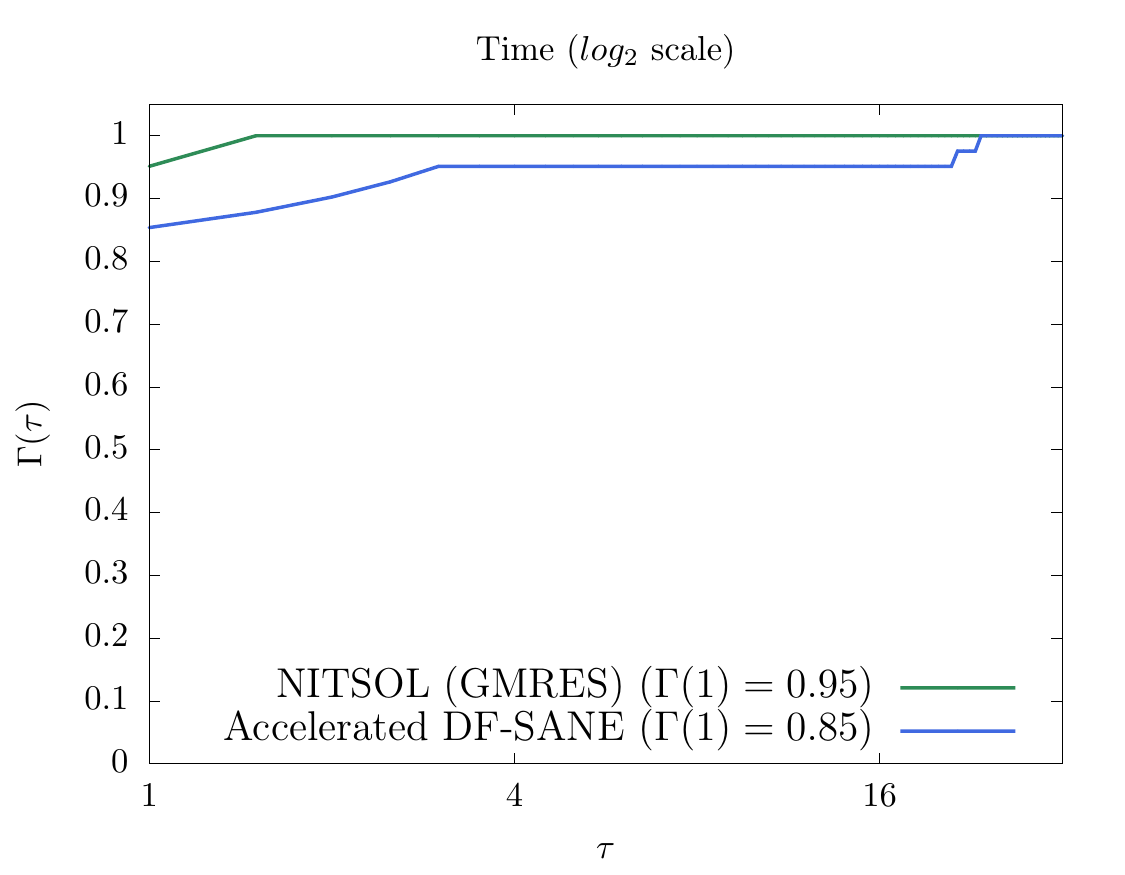}
\caption{Performance profiles of Accelerated DF-SANE (in \proglang{Fortran}) and \pkg{NITSOL} (GMRES) considering the~41 problems from the \pkg{CUTEst} collection in which both methods found a solution.}
\label{ppdiisnitsol}
\end{figure}

In a second experiment, in order to put our method in perspective relatively to a method that represents the state of the art in solving nonlinear systems, we compared Accelerated DF-SANE with \pkg{NITSOL} (GMRES). Since \pkg{NITSOL} (GMRES) is coded in \proglang{Fortran}, we considered the \proglang{Fortran} version of Accelerated DF-SANE in this comparison. Of course, we considered \pkg{NITSOL} (GMRES) \textit{without} Jacobians. Table~\ref{accvsnitsol} and Figure~\ref{ppdiisnitsol} show the results. As in Table~\ref{accvsravi}, in column $\|F(x_*)\|_2$, figures in red are the ones that \textit{do not} satisfy~\eqref{scsuccess}. In all cases the final iterate of \pkg{NITSOL} (GMRES) does not satisfy~\eqref{scsuccess}, \pkg{NITSOL} (GMRES) stops by ``too small step in a line search'' (flag equal to~6). 

Figures in Table~\ref{accvsnitsol} show that both Accelerated DF-SANE and \pkg{NITSOL} (GMRES) solve~45 problems. There are~41 problems that were solved by both methods, 4 problems that were solved by Accelerated DF-SANE only, and~4 problems that were solved by \pkg{NITSOL} (GMRES) only.  So, both methods appear to be equally robust. 

As well as Figure~\ref{ppdiisravi}, Figure~\ref{ppdiisnitsol} focuses on efficiency and, thus, it considers only the 41 problems in which both, Accelerated DF-SANE and \pkg{NITSOL} (GMRES), found a solution. Figure~\ref{ppdiisnitsol}(a) considers number of functional evaluations as performance metric; while Figure~\ref{ppdiisnitsol}(b) considers CPU time. Figure~\ref{ppdiisnitsol}(a) shows that \pkg{NITSOL} (GMRES) used less functional evaluations in 63\% of the problems; while Accelerated DF-SANE used less functional evaluations in 39\% of the problems. (The sum of the percentages is  slightly larger than 100\% because ties are counted twice.) The fact that the two curves reach~0.9 before $\tau=10$ means that in~90\% of the problems the number of function evaluations is of the same order. The Accelerated DF-SANE curve reaches the value 1 for $\tau>1000$ due to only~3 problems. In the problems \textsc{Recipe}, \textsc{Heart8} and \textsc{Hatfldg}, Accelerated DF-SANE consumes approximately 14, 33 and 1790 times more function evaluations than \pkg{NITSOL} (GMRES). On the other hand, the curve of \pkg{NITSOL} (GMRES) reaches the value 1 between $\tau=10$ and $\tau=100$ because in the problem \textsc{Waysea1ne}  \pkg{NITSOL} (GMRES) uses~41 times more function evaluations than Accelerated DF-SANE.

The performance profile of the Figure~\ref{ppdiisnitsol}(b) that considers CPU time as performance measure, shows a similar scenario, contaminated by the fact of having a large proportion of small problems. The figure says that \pkg{NITSOL} (GMRES) is the fastest method in 95\% of the problems; while Accelerated DF-SANE is the fastest method in 85\% of the problems, i.e., there are a lot of ties. (As it can be observed in Table~\ref{accvsnitsol}, approximately 90\% of the CPU times associated with problems that are solved by both methods are smaller than 0.1 seconds; and CPU times smaller than 0.01 seconds are considered ties.) The curve of \pkg{NITSOL} (GMRES) reaches~1 before $\tau=2$ because in no problem does \pkg{NITSOL} (GMRES) uses more than twice the time of Accelerated DF-SANE. Accelerated DF-SANE also did not use more than twice the time of NITSOL in~37 out of the~41 problems. On the remaining~4 problems, Accelerated DF-SANE uses a little more than twice as much time on \textsc{Chandheu} and \textsc{Oscigrne} (which is why the curve passes~0.95 before $\tau=3$) and on problems \textsc{Heart8} and \textsc{Hatfldg} it uses~21 and~23 times as much time. 

Summing up, we conclude that, while both methods are equally robust, \pkg{NITSOL} (GMRES) is slightly more efficient that Accelerated DF-SANE in the considered set of problems. On the other hand, it is worth noticing that numerical experiments in~\cite{bmdiis} showed that Accelerated DF-SANE outperforms \pkg{NITSOL} (GMRES) to a large extent on an important class of large-scale problems coming from the discretization of partial differential equations. Of course, the opposite situation can also occur, which justifies the availability of both methods.

\begin{center}
{\footnotesize\tabcolsep=2pt
\begin{longtable}{lr@{\hskip 2ex}rrrrrrrrr}
\toprule
\multirow{2}{*}{Problem} &
\multirow{2}{*}{$n$} &
\multicolumn{4}{c}{Accelerated DF-SANE} & &
\multicolumn{4}{c}{\pkg{NITSOL} (GMRES)}\\ 
\cmidrule{3-6} \cmidrule{8-11}
&        & $\|F(x_*)\|$&  \# iter  &  \# feval  &     time  && $\|F(x_*)\|$&  \# iter & \# feval &      time  \\ 
\midrule
     BOOTH &      2 &       $9.9$E$-16$ &          2 &          7 &   0.000014 &&       $4.6$E$-09$ &          3 &          8 &   0.000039 \\  
   CLUSTER &      2 &       $8.3$E$-07$ &         23 &        108 &   0.000048 &&       $1.2$E$-09$ &          9 &         25 &   0.000046 \\  
    CUBENE &      2 &       $4.0$E$-13$ &          9 &         20 &   0.000022 &&       $2.1$E$-10$ &         38 &        108 &   0.000076 \\  
DENSCHNCNE &      2 &       $2.3$E$-11$ &         10 &         23 &   0.000029 &&       $6.7$E$-07$ &          6 &         15 &   0.000043 \\  
DENSCHNFNE &      2 &       $2.7$E$-07$ &          7 &         23 &   0.000019 &&       $1.6$E$-13$ &          5 &         16 &   0.000044 \\  
  FREURONE &      2 &       $1.5$E$-08$ &         16 &         55 &   0.000025 && \red{$7.0$E$+00$} &         16 &        112 &   0.000058 \\  
    GOTTFR &      2 &       $1.3$E$-07$ &         23 &         67 &   0.000031 &&       $3.6$E$-09$ &         70 &        236 &   0.000133 \\  
  HIMMELBA &      2 &       $0.0$E$+00$ &          2 &          7 &   0.000013 &&       $2.5$E$-08$ &          3 &          8 &   0.000042 \\  
  HIMMELBC &      2 &       $8.4$E$-08$ &          5 &         13 &   0.000018 &&       $1.1$E$-06$ &          6 &         14 &   0.000041 \\  
  HIMMELBD &      2 & \red{$2.4$E$+00$} &   11577102 &  439522728 & 180.000000 && \red{$2.4$E$+00$} &         48 &        246 &   0.000164 \\  
       HS8 &      2 &       $4.4$E$-08$ &          5 &         13 &   0.000018 &&       $2.4$E$-11$ &         11 &         24 &   0.000045 \\  
    HYPCIR &      2 &       $8.7$E$-10$ &          6 &         14 &   0.000017 &&       $5.2$E$-07$ &          5 &         13 &   0.000041 \\  
  POWELLBS &      2 &       $1.4$E$-06$ &   54229896 & 1259707609 & 152.775132 && \red{$1.9$E$-06$} &        231 &        692 &   0.000206 \\  
  POWELLSQ &      2 & \red{$1.4$E$-00$} &   13690098 &   34211713 & 180.000000 && \red{$1.3$E$+00$} &      37498 &     309809 &   0.044447 \\  
  PRICE3NE &      2 &       $3.9$E$-10$ &          7 &         19 &   0.000020 &&       $4.4$E$-10$ &          7 &         20 &   0.000046 \\  
  PRICE4NE &      2 &       $1.3$E$-10$ &         10 &         27 &   0.000030 &&       $3.0$E$-09$ &         10 &         27 &   0.000048 \\  
   RSNBRNE &      2 &       $2.2$E$-16$ &         56 &        204 &   0.000054 &&       $1.4$E$-06$ &         55 &        161 &   0.000075 \\  
  SINVALNE &      2 &       $4.9$E$-15$ &         16 &         77 &   0.000040 &&       $1.9$E$-14$ &          6 &         19 &   0.000042 \\  
 WAYSEA1NE &      2 &       $1.3$E$-10$ &         12 &         36 &   0.000023 &&       $3.4$E$-08$ &        331 &       1485 &   0.000291 \\  
 WAYSEA2NE &      2 &       $8.4$E$-07$ &        481 &       2179 &   0.000401 &&       $1.3$E$-09$ &        766 &       3751 &   0.000677 \\  
DENSCHNDNE &      3 &       $2.3$E$-07$ &         26 &         62 &   0.000043 &&       $1.5$E$-06$ &         22 &         71 &   0.000065 \\  
DENSCHNENE &      3 &       $9.6$E$-11$ &          6 &         16 &   0.000032 &&       $1.5$E$-09$ &          7 &         19 &   0.000046 \\  
   HATFLDF &      3 &       $1.4$E$-08$ &         26 &         78 &   0.000049 &&       $9.6$E$-07$ &         71 &        233 &   0.000117 \\  
HATFLDFLNE &      3 & \red{$7.9$E$-03$} &   11587628 &  252488903 & 180.000000 && \red{$7.8$E$-03$} &        372 &       2843 &   0.000672 \\  
   HELIXNE &      3 &       $2.8$E$-09$ &         13 &         35 &   0.000045 && \red{$5.0$E$+01$} &          0 &         14 &   0.000040 \\  
  HIMMELBE &      3 &       $9.7$E$-16$ &          9 &         21 &   0.000023 &&       $7.3$E$-09$ &          2 &          9 &   0.000043 \\  
    RECIPE &      3 &       $6.2$E$-07$ &         72 &        403 &   0.000116 &&       $1.4$E$-06$ &         10 &         28 &   0.000048 \\  
  ZANGWIL3 &      3 &       $1.4$E$-14$ &          3 &         11 &   0.000015 &&       $5.2$E$-07$ &          3 &         10 &   0.000045 \\  
  POWELLSE &      4 &       $7.3$E$-07$ &         24 &         70 &   0.000061 &&       $1.5$E$-06$ &         13 &         61 &   0.000064 \\  
POWERSUMNE &      4 & \red{$1.2$E$-02$} &    8695243 &  130633973 & 180.000000 &&       $1.6$E$-06$ &       1417 &       7084 &   0.004665 \\  
    HEART6 &      6 &       $1.5$E$-06$ &     124382 &    1818751 &   0.561869 && \red{$2.7$E$-01$} &       3854 &      28811 &   0.013372 \\  
    HEART8 &      8 &       $2.8$E$-06$ &     181971 &    2866905 &   0.993949 &&       $2.2$E$-06$ &      11360 &      86495 &   0.046512 \\  
  COOLHANS &      9 &       $1.5$E$-06$ &         10 &         45 &   0.000056 &&       $2.3$E$-06$ &          7 &         22 &   0.000057 \\  
  MOREBVNE &     10 &       $1.6$E$-06$ &         37 &        219 &   0.000124 &&       $7.9$E$-08$ &          4 &         33 &   0.000068 \\  
  OSCIPANE &     10 & \red{$1.0$E$+00$} &    8608149 &  322784536 & 180.000000 && \red{$1.0$E$+00$} &       2411 &      50650 &   0.022718 \\  
 TRIGON1NE &     10 &       $1.9$E$-06$ &         13 &         29 &   0.000063 &&       $2.5$E$-06$ &          5 &         26 &   0.000069 \\  
   INTEQNE &     12 &       $9.2$E$-07$ &          3 &          7 &   0.000021 &&       $3.3$E$-07$ &          4 &         10 &   0.000065 \\  
   HATFLDG &     25 &       $5.0$E$-06$ &      22708 &     356246 &   0.232828 &&       $7.8$E$-07$ &         44 &        199 &   0.000263 \\  
   HYDCAR6 &     29 & \red{$5.0$E$-03$} &    2661134 &   61551663 & 180.000000 && \red{$3.3$E$-01$} &         30 &        781 &   0.002596 \\  
  METHANB8 &     31 & \red{$1.2$E$-04$} &    2577703 &   67500153 & 180.000000 && \red{$1.4$E$-02$} &          6 &        472 &   0.001542 \\  
  METHANL8 &     31 & \red{$4.4$E$-03$} &    2764968 &   66380772 & 180.000000 && \red{$6.1$E$-01$} &         28 &       1052 &   0.003356 \\  
  HYDCAR20 &     99 & \red{$3.9$E$-02$} &     917448 &   19172981 & 180.000000 && \red{$9.2$E$+00$} &          3 &        287 &   0.003190 \\  
  LUKSAN21 &    100 &       $8.9$E$-06$ &         48 &        441 &   0.001177 &&       $6.1$E$-06$ &         17 &        123 &   0.000562 \\  
 MANCINONE &    100 &       $5.9$E$-07$ &          5 &         17 &   0.009272 &&       $3.9$E$-06$ &          4 &         11 &   0.005929 \\  
    QINGNE &    100 &       $4.8$E$-06$ &         21 &         45 &   0.000233 &&       $4.3$E$-06$ &         10 &         35 &   0.000150 \\  
   ARGTRIG &    200 &       $1.2$E$-05$ &         57 &        199 &   0.016417 &&       $1.1$E$-05$ &          5 &         86 &   0.007244 \\  
  BROWNALE &    200 &       $1.0$E$-05$ &          9 &         25 &   0.001325 &&       $3.1$E$-07$ &          3 &          9 &   0.000512 \\  
  CHANDHEU &    500 &       $1.4$E$-05$ &         18 &         99 &   0.140877 &&       $1.5$E$-05$ &         10 &         51 &   0.065100 \\  
  10FOLDTR &   1000 & \red{$2.2$E$+07$} &       9445 &     272830 & 180.000000 &&       $2.7$E$-05$ &         54 &       6563 &   4.562871 \\  
       KSS &   1000 &       $9.3$E$-06$ &          5 &         17 &   0.023044 &&       $2.2$E$-08$ &          6 &         13 &   0.017676 \\  
    MSQRTA &   1024 & \red{$4.7$E$+01$} &      68938 &    1137480 & 180.000000 && \red{$5.5$E$+01$} &         17 &       1351 &   0.210034 \\  
    MSQRTB &   1024 & \red{$4.6$E$+01$} &      61153 &    1138024 & 180.000000 && \red{$5.9$E$+01$} &         13 &       1964 &   0.306907 \\  
   EIGENAU &   2550 & \red{$1.6$E$+02$} &      12625 &     234607 & 180.000000 && \red{$1.6$E$+02$} &         17 &        850 &   0.768981 \\  
    EIGENB &   2550 & \red{$9.6$E$+00$} &      15297 &     234454 & 180.000000 && \red{$9.8$E$+00$} &          9 &        382 &   0.361665 \\  
    EIGENC &   2652 & \red{$9.2$E$+01$} &      14864 &     218919 & 180.000000 && \red{$9.7$E$+01$} &         33 &       2169 &   2.097641 \\  
NONMSQRTNE &   4900 & \red{$2.4$E$+02$} &       5731 &      85005 & 180.000000 && \red{$2.3$E$+02$} &         23 &        915 &   1.804071 \\  
  BROYDN3D &   5000 &       $5.3$E$-05$ &         12 &         25 &   0.005502 &&       $2.8$E$-05$ &          5 &         19 &   0.002987 \\  
  BROYDNBD &   5000 & \red{$2.4$E$+00$} &      58861 &     934685 & 180.000000 && \red{$7.7$E$+00$} &         11 &        607 &   0.176834 \\  
  BRYBNDNE &   5000 & \red{$2.4$E$+00$} &      57595 &     915686 & 180.000000 && \red{$7.7$E$+00$} &         11 &        607 &   0.176482 \\  
  NONDIANE &   5000 & \red{$1.0$E$+00$} &      83049 &    1603628 & 180.000000 && \red{$6.1$E$+02$} &        686 &      10094 &   1.873028 \\  
 SBRYBNDNE &   5000 & \red{$2.5$E$+02$} &      45364 &     906538 & 180.000000 && \red{$2.7$E$+02$} &         50 &       2935 &   0.918074 \\  
SROSENBRNE &   5000 &       $2.5$E$-09$ &          9 &         34 &   0.004332 &&       $2.1$E$-08$ &          4 &         11 &   0.001462 \\  
SSBRYBNDNE &   5000 & \red{$1.7$E$+02$} &      50681 &     944507 & 180.000000 && \red{$1.6$E$+02$} &        128 &       9043 &   2.794424 \\  
TQUARTICNE &   5000 & \red{$8.3$E$-01$} &     175237 &    1886434 & 180.000000 &&       $1.5$E$-07$ &          2 &          6 &   0.000899 \\  
  OSCIGRNE & 100000 &       $1.8$E$-04$ &         28 &         66 &   0.461298 &&       $1.5$E$-04$ &          7 &         34 &   0.158588 \\  
   CYCLIC3 & 100002 & \red{$6.2$E$-01$} &       3011 &      53186 & 180.000000 &&       $1.7$E$-04$ &        282 &        992 &   4.070610 \\  
  YATP1CNE & 123200 &       $2.6$E$-07$ &         14 &         41 &   0.889454 &&       $1.4$E$-04$ &         17 &         48 &   0.970848 \\  
   YATP1NE & 123200 &       $2.6$E$-07$ &         14 &         41 &   0.891586 &&       $1.4$E$-04$ &         17 &         48 &   0.974741 \\  
  YATP2CNE & 123200 & \red{$3.1$E$+04$} &        800 &      12314 & 180.000000 &&          \red{--} &         -- &         -- & 180.000000 \\  
   YATP2SQ & 123200 & \red{$4.1$E$+04$} &        791 &      12362 & 180.000000 &&          \red{--} &         -- &         -- & 180.000000 \\  
\bottomrule
\caption{Detailed results of the application of Accelerated DF-SANE (in \proglang{Fortran}) and \pkg{NITSOL} (GMRES) to the~70 considered problems from the \pkg{CUTEst} collection.}
\label{accvsnitsol}
\end{longtable}
}
\end{center}
\vspace{-3em}

A side note comparing the \proglang{R} and \proglang{Fortran} implementations of Accelerated DF-SANE is in order. Comparing Tables~\ref{accvsravi} and~\ref{accvsnitsol}, it can be seen that they deliver slightly different results in a few problems and deliver identical results in~40 problems out of the~44 problems in which none of the versions stops by reaching the CPU time limit. If we consider these~40 problems, in which both versions performed an identical number of iterations and functional evaluations, the \proglang{Fortran} version uses, in average, around 10\% of the CPU time required by the \proglang{R} version of the method.

\section{Conclusions} \label{conclusions}

In~\cite{bmdiis}, where it was shown that an acceleration scheme based on the Sequential Secant Method could improve the performance of the derivative-free spectral residual method~\citep{lmr}, numerical experiments with very large problems coming from the discretization of partial differential equations were presented. In the considered family of problems, Accelerated DF-SANE outperformed DF-SANE and \pkg{NITSOL} (GMRES) by a large extent. 

In the present work, an \proglang{R} implementation of the method proposed in~\cite{bmdiis} was introduced. In addition, numerical experiments considering \textit{all} nonlinear systems of equations from the well-known \pkg{CUTEst} collection were presented. Default dimensions of the problems were considered; and the collection includes small-, medium-, and large-scale problems. Results shown that the proposed method is much more robust than the DF-SANE method included in the \proglang{R} package \pkg{BB}~\citep{ravi}; while it is as robust and almost as efficient as the state-of-the-art classical \pkg{NITSOL} (GMRES) method (coded in \proglang{Fortran}). Therefore, the proposed method appears as a useful and robust alternative for solving nonlinear systems of equations without derivatives to the users of \proglang{\proglang{R}} language. 

As a byproduct, an interface to test derivative-free nonlinear systems solvers developed in~\proglang{R} with the widely-used test problems from the \pkg{CUTEst} collection~\citep{cutest} was also provided.

\bibliographystyle{plain}
\bibliography{bgmmdiis}

\end{document}